\documentclass{elsarticle}

\usepackage{graphics}

\usepackage{amssymb}

\usepackage[nodots]{numcompress}
\usepackage{psfrag}
\usepackage{color}
\usepackage{verbatim}
\usepackage{import}
\usepackage{amsmath,amsthm,mathrsfs,epsf,a4,color,graphicx,eucal}

\usepackage{subfigure}
\newcommand{\COL}[1]{{#1}}
\newcommand{\COLL}[1]{{#1}}

\newcommand{\REM}[1]{}

\newcommand\R{\mathbb R}
\newcommand\N{\mathbb N}

\newcommand{\eq}[1]{(\ref{#1})}
\newcommand\DT[1]{\mathchoice
                 {{\buildrel{\hspace*{.1em}\text{\LARGE.}}\over{#1}}}
                 {{\buildrel{\hspace*{.1em}\text{\Large.}}\over{#1}}}
                 {{\buildrel{\hspace*{.1em}\text{\large.}}\over{#1}}}
                 {{\buildrel{\hspace*{.1em}\text{\large.}}\over{#1}}}}
\newcommand\DDT[1]{\mathchoice
   {{\buildrel{\hspace*{.1em}\text{\LARGE.\hspace*{-.1em}.}}\over{#1}}}
   {{\buildrel{\hspace*{.1em}\text{\Large.\hspace*{-.1em}.}}\over{#1}}}
   {{\buildrel{\hspace*{.1em}\text{\large.\hspace*{-.1em}.}}\over{#1}}}
   {{\buildrel{\hspace*{.1em}\text{\large.\hspace*{-.1em}.}}\over{#1}}}}

\newcommand\bbC{\mathbb C}
\newcommand\bbD{\mathbb D}

\renewcommand\d{\mathrm d}

\newcommand\calE{\mathscr E}
\newcommand\calG{\mathscr G}

\newcommand\calR{\mathscr R}

\newcommand\GC{\Gamma_{\!\mbox{\tiny\rm C}}}

\newcommand{\GDir}{\Gamma_{\!\mbox{\tiny\rm D}}}
\newcommand{\GNeu}{\Gamma_{\!\mbox{\tiny\rm N}}}

\newcommand{\Gdir}{\Gamma_{\!\mbox{\tiny\rm D}}}
\newcommand{\Gnew}{\Gamma_{\!\mbox{\tiny\rm N}}}

\newcommand{\dd}{\,\mathrm{d}}

\journal{International Journal of Solids and Structures}

\begin{document}

\begin{frontmatter}


\title{A simple and efficient BEM implementation of quasistatic linear visco-elasticity}
\author{C.G.~Panagiotopoulos}

\author{V.~Manti\v{c}\corref{cor1}}
\ead{mantic@us.es}
\cortext[cor1]{Corresponding author. Tel.:+34-954-482135; fax:+34-954-461637}
\address{Group of Elasticity and Strength of Materials, Department of Continuum Mechanics \\ School of Engineering,
University of Seville\\[-0.25em]  Camino de los Descubrimientos s/n, ES-41092 Sevilla, Spain}

\author{T.Roub\'\i\v cek}
\address{Mathematical Institute, Charles University, Sokolovsk\'a 83,
CZ--18675 Praha~8, Czech Republic\\
Institute of Thermomechanics of the ASCR, Dolej\v skova~5,
CZ--18200 Praha 8, Czech Republic}


\address{}

\begin{abstract}
A simple, yet efficient procedure to solve quasistatic problems of
special linear visco-elastic solids  at small strains with equal
rheological response in all tensorial components, utilizing boundary element method (BEM), is introduced.
This procedure is based on the implicit discretisation in time (the so-called Rothe method) combined with  a simple ``algebraic'' transformation of variables, leading to a numerically stable procedure (proved explicitly by discrete energy estimates), which can be easily implemented in a  BEM code
to solve initial-boundary value visco-elastic problems  by using the Kelvin
elastostatic fundamental solution only. It is worth mentioning that no   inverse Laplace transform is required here. The formulation is straightforward
for both 2D and 3D problems  involving  unilateral
frictionless contact. Although the focus is to the simplest  Kelvin-Voigt
rheology, a generalization to Maxwell, Boltzmann, Jeffreys, and Burgers
rheologies is proposed, discussed, and implemented in the BEM code too.
A few  2D and 3D initial-boundary value problems, one of them  with unilateral
frictionless contact, are solved numerically.

\end{abstract}

\begin{keyword}
     boundary element method
\sep implicit time discretisation
\sep quasistatic linear visco-elasticity
\sep unilateral contact
\sep Kelvin-Voigt rheology
\sep Maxwell rheology
\sep standard linear solids
\sep Jeffreys rheology
\sep Burgers rheology



\end{keyword}

\end{frontmatter}


\section{Introduction}\label{sec:Intro}
A large number of engineering  and (e.g.\ geo-)physical applications
consider materials that
 exhibit
{\it visco-elastic behaviour}. A typical example of such a behaviour is the
mechanical response of polymers and polymer-matrix composites, or rocks
undergoing aseismic slip, etc.
Visco-elasticity accounts for the dependence of stresses and strains on time,
and response of real visco-elastic solids or structures is usually analysed
numerically by the finite    or   boundary element
methods (FEM or BEM). When inertial effects are neglected, usually because of sufficiently slow
external loading, the
model is addressed as {\it quasistatic}.  The quasistatic linear
visco-elasticity theory provides a usable engineering approximation for
many applications in polymer and composites engineering, among others. There
are several   models  describing   visco-elastic behaviour of materials
obtained by a generalization of simple 1D models to 2D or
3D ones. One of these well-known models, often adopted in
designing procedures,  is the Kelvin-Voigt model.

There are   four main approaches  to  {\it quasistatic linear visco-elastic analysis by BEM}. The first  and most commonly applied  approach uses the correspondence principle to establish an associated elastic problem solved in the Laplace transform domain. Then, the solution in  time domain is recovered by a numerical inversion~\cite{Rizzo,Manolis_Beskos,Kusama_Mitsui,Sladeks1984,Carini1986,ChenHwu2011}. \COL{The second approach works directly in the time domain, however, it requires a time dependent fundamental solution}~\cite{LeeWestmann,Cezario2011,Zhu2011}.  The third, a kind of mixed, approach  also solves the problem in time domain, but  uses  the Laplace transformed fundamental solutions with a convolution quadrature leading to a time stepping procedure  without the knowledge of the time dependent fundamental solution~\cite{Schanz, SchanzAntes, Syngellakis}. The fourth, a kind of direct, approach which utilizes the Kelvin  elastostatic fundamental solution was introduced by Mesquita and Coda in \cite{MesquitaCoda2001,MesCod02} for both Kelvin-Voigt
and Boltzmann visco-elastic models. The Somigliana displacement and stress
indentities are rewritten  to obtain \emph{visco-elastic boundary-integral-representations} (BIRs) for these models. After the BEM discretisation of these BIRs, a
{\it finite difference approximation} of velocities leads to a time marching
scheme. This approach was later applied to the problem of circular holes and
elastic inclusions in a visco-elastic plane \cite{Huang2005EABE,Huang2005CM}. A   brief
presentation of several BEM procedures for problems of visco-elasticity may
  be found in \cite{Marques_Creus}.

The novelty of the present approach consists in a particular application of
the Rothe method (i.e.\ the time discretisation by the implicit Euler
formula, cf.\ e.g.\ \cite{Roub13NPDE}) to the governing partial differential
equations (PDE), where after this
time discretisation,
a suitable variable transform is carried out to convert it in each time
step to a linear {\it auxiliary
elastostatic problem}
\REM{(Tom: I call it ``auxiliary elastostatic'', OK?)}
with proper boundary conditions. Once this linear
elastostatic problem is solved the actual displacements, stresses and strains of the visco-elastic problem in this time step are recovered and used in the next step, an efficient recursive procedure being obtained in this way. For the sake of simplicity of explanation, the main steps of the procedure proposed are first explained for the simple Kelvin-Voigt model, and then briefly generalized to other basic linear visco-elastic rheologies. The present procedure can be implemented in any elastostatic FEM or BEM code. The present work is based on the collocation BEM formulation due to its advantages as no domain variables appear in the \COL{solution of the problem}. Additionally the stability of the present time
discretisation can be established.
Although there are evident similarities with the previous work by Mesquita and Coda, the present theoretical formulation is much straightforward  showing  in  a more transparent way that any linear elastostatic BEM code  can be applied  to linear visco-elastic analysis requiring just minor modifications.

Under these assumptions, the purpose of this work is to present and numerically verify  a simple yet efficient methodology  for BEM analysis of quasistatic visco-elastic solids, initially scrutinizing the Kelvin-Voigt material
in Sections~\ref{sect-KV}-\ref{sec:BEF} and later, in
Section~\ref{sect-other}, further extended to other models
usually found in engineering or physical applications. The approach may
be considered as
a time domain one, where no special time-depended fundamental solution, neither domain integration, is needed. Another important engineering problem  treated in this work  is a contact of visco-elastic bodies \cite{Graham}.
\section{The mixed unilateral initial-boundary-value problem for Kelvin-Voigt visco-elastic body}\label{sect-KV}

\begin{figure}[h]
   \centering
   \def\svgwidth{0.5\columnwidth}
   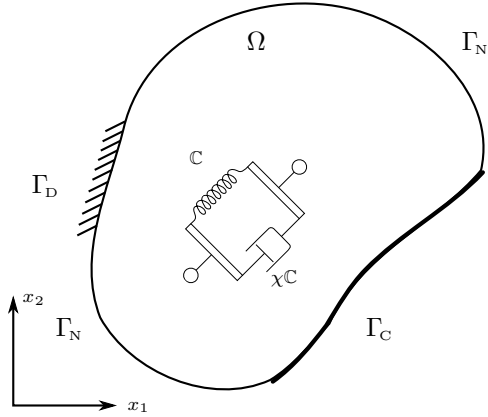
   \caption{2D schematic illustration of the geometry and
notation of the boundary-value problems considered. In the bulk, a
visco-elastic rheology from Fig.~\ref{fig-KVMA} is schematically depicted.}
   \label{fig_prob_conf}
\end{figure}

The following boundary-value problem on a
domain $\Omega\subset\R^d$, $d=2,3$,
is used in the subsequent developments, where also the standard model of the frictionless unilateral
Signorini contact is considered, see
Figure~\ref{fig_prob_conf},
\begin{subequations}\label{statement}
\begin{align}
\label{statement_01}
&&&\mathrm{div}\,\bbC\epsilon+f=0\ \ \ \ \text{ with }\ \
\epsilon=\epsilon(u,\DT{u})=e(u{+}\chi\DT{u})\qquad &&\text{on }\Omega,&&&&\\
\label{statement_02}
&&&u=w &&\text{on }\GDir,\\
\label{statement_03}
&&&\mathfrak{t}(\epsilon)=\big(\bbC\epsilon\big)\big|_{\Gamma}\vec{n}=g &&\text{on }\GNeu,\\
\label{statement_04}
&&&u{\cdot}\vec{n}\le 0,\ \ \ \ \mathfrak{t}_{\rm n}(\epsilon)\le0,\ \ \ \
(u{\cdot}\vec{n})\mathfrak{t}_{\rm n}(\epsilon)=0,\ \ \ \ \mathfrak{t}_{\rm t}(\epsilon)=0
&&\text{on }\GC,
\end{align}
where $u$ is the displacement and
$e=e(u)=\frac12(\nabla u)^\top\!{+}\frac12\nabla u$ the small-strain tensor,
and $\bbC$ is the fourth order tensor of elastic moduli,
while $\chi>0$ a given relaxation time. Furthermore,
$\vec{n}=\vec{n}(\vec{x})$ is the unit outward normal to
$\Gamma=\partial\Omega$ at $x$, $\mathfrak{t}_{\rm n}(\epsilon)=\mathfrak{t}(\epsilon){\cdot}\vec{n}$,
and $\mathfrak{t}_{\rm t}(\epsilon)=\mathfrak{t}(\epsilon)-\mathfrak{t}_{\rm n}(\epsilon)\vec{n}$. It is straightforward to generalize the above problem formulation and  all the results below to  several (visco-)elastic solids in contact with a non-negative gap defined at a possible contact zone  $\GC$ (see Example~\ref{Contact}). Actually, pertinent indications in this sense will be given at some places below.
We further consider the initial-value problem
for (\ref{statement}a-d) for time $t>0$  by prescribing the initial condition
\begin{align}
\label{IC}
u(0)=u_0.
\end{align}
\end{subequations}

\begin{figure}[h]
   \centering
   \vspace{2em}
   \def\svgwidth{0.25\columnwidth}
   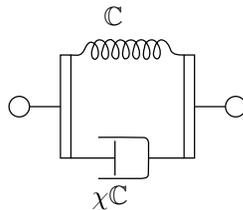
   \caption{Mechanical analog of Kelvin-Voigt model.}
   \label{fig-KVMA}
\end{figure}

The mechanical 1D analog of the above model is shown in Figure~\ref{fig-KVMA}. According to this figure, since the two components of the model are arranged in parallel, the strains in each component are identical and equal to $e(u)$, while for the stress it holds,
\begin{align}\label{Kelvin}
\sigma = \bbC e(u)+\chi\bbC e(\DT{u}),
\end{align}
where the actual or total stress field is defined as the sum of the elastic
and visco-elastic part.
 The  Kelvin-Voigt model is known to be very
 effective for predicting
creep, but less at describing the relaxation behavior.
For this reason other advanced and more complex rheological models
exploiting auxiliary internal parameters have been defined and used.
Eliminating these internal parameters leads to higher order time derivatives
involved in the model, cf.\ Section~\ref{sect-other}.

\section{Discretisation in time and space
}\label{sec:BEF}

\COL{
We perform the discretisation of the initial-boundary value problem
\eqref{statement} by the implicit formula in time and by
the boundary-element method in space.}

\subsection{Time discretisation}
Using an equidistant partition of the time
interval $[0,T]$ with a time step $\tau>0$ such that $T/\tau\in\N$, we
consider:
\begin{subequations}\label{timedisc}
\begin{align}
\label{timedisc_01}
&&&\mathrm{div}\,\bbC\epsilon_\tau^k+f_\tau^k=0\quad
\text{ with }\ \epsilon_\tau^k=
e\big(u_\tau^k+\chi(u_\tau^k{-}u_\tau^{k-1})/\tau\big)
 &&\text{on }\Omega,&&&&\\
\label{timedisc_02}
&&&u_\tau^k=w_\tau^k &&\text{on }\GDir,\\
\label{timedisc_03}
&&&\mathfrak{t}(\epsilon_\tau^k)=\big(\bbC\epsilon_\tau^k\big)\big|_{\Gamma}\vec{n}=g_\tau^k &&\text{on }\GNeu,\\
\label{timedisc_04}
&&&u_\tau^k{\cdot}\vec{n} \le 0,\ \ \ \
\mathfrak{t}_{\rm n}(\epsilon_\tau^k)\le0,\ \ \ \
(u_\tau^k{\cdot}\vec{n})\mathfrak{t}_{\rm n}(\epsilon_\tau^k)=0,\ \ \ \ \mathfrak{t}_{\rm t}(\epsilon_\tau^k)=0
&&\text{on }\GC,
\end{align}
with $w_\tau^k=w(k\tau)$, $f_\tau^k=f(k\tau)$ and $g_\tau^k=g(k\tau)$, and
proceed recursively for $k=1,...,T/\tau$ with starting for $k=1$ from
\begin{align}
u_\tau^0=u_0.
\end{align}
\end{subequations}

This implicit time discretisation is numerically stable in the
sense that the discrete solution $u_\tau^k$ stays bounded if $\tau\to0$ in
a suitable norm provided the data $u_0$, $f$, and $g$ are qualified
appropriately. More specifically, this can be seen from the
discrete variant \COL{(as an upper inequality) of the continuous
energy-conservation equality   \eqref{TotalEngr}, introduced and discussed in Appendix,}
i.e.
\begin{align}\label{TotalEngr+disc}
&\calE(u_\tau^k)+\sum_{l=1}^k
\int_{\Omega}\!\chi\bbC e\Big(\frac{u_\tau^l{-}u_\tau^{l-1}}\tau\Big){:}
e\Big(\frac{u_\tau^l{-}u_\tau^{\COLL{l}-1}}\tau\Big)\,\d x
\nonumber\\&\qquad
\le\calE(u_0)+\sum_{l=1}^k\bigg(\int_{\Omega}\!
f_\tau^l{\cdot}\frac{u_\tau^l{-}u_\tau^{l-1}}\tau\,\d x+
\int_{\GNeu}\!\!g_\tau^l{\cdot}\frac{u_\tau^l{-}u_\tau^{l-1}}\tau\dd S\bigg).
\end{align}
The inequality in \eqref{TotalEngr+disc}
rely on convexity of the stored energy $\calE$.
%


\subsection{Transform of the visco-elastic to an auxiliary elastic-like problem}
BEM standardly uses the so-called boundary integral
operators which
are explicitly known in specific static cases, here for the homogeneous linear elastic 
material which we consider in what follows. Yet, we have to calculate
visco-elastic modification and here we benefit from choosing
the ansatz of the tensor of viscous moduli as simply proportional
to the elastic moduli, i.e.\ $\chi\bbC$. Therefore
we can use BEM with the same
boundary integral operators as in the static case utilizing a transformation originally proposed in \cite{Ro_AdhVisc} and numerically implemented in \cite{RoPaMa_Visc},
by defining a new  auxiliary variable, in view of \eqref{statement_01}, as 
\begin{align}\label{eqauxfield}
v_{\tau}^{k}=u_\tau^k+\chi\frac{u_\tau^k{-}u_\tau^{k-1}}\tau.
\end{align}
In terms of this new variable, one obviously has
the Kelvin-Voigt strain
$\epsilon_\tau^k=e(v_\tau^k)$, the velocity $(u_\tau^k{-}u_\tau^{k-1})/\tau=
(v_\tau^k{-}u_\tau^{k-1})/(\tau{+}\chi)$, and the displacement recovered by
\begin{align}\label{eqphysicalfield}
u_\tau^k=(\tau v_\tau^k{+}\chi u_\tau^{k-1})/(\tau{+}\chi),
\end{align}
which is to be used in \eqref{timedisc_01}-\eqref{timedisc_03}, where we assume $\GC=\emptyset$, leading to the transformed time discretized problem
\begin{subequations}\label{time-space-disc}
\begin{align}
\label{time-space-disc_01}
&&&\mathrm{div}\,\bbC e(v_\tau^k)+f_\tau^k=0\qquad &&\text{on }\Omega,&&&&\\
\label{time-space-disc_02}
&&&v_\tau^k=\frac{\chi{+}\tau}\tau w_\tau^k-\frac\chi\tau w_\tau^{k-1}  &&\text{on }\GDir,\\
\label{time-space-disc_03}
&&&\mathfrak{t}(e(v_\tau^k))=\big(\bbC e(v_\tau^k)\big)\big|_{\Gamma}\vec{n}=g_\tau^k &&\text{on }\GNeu,
\end{align}
\end{subequations}
with $u_\tau^{k-1}=(\tau v_\tau^{k-1}{+}\chi u_\tau^{k-2})/(\tau{+}\chi)$,
and proceeding recursively for $k=1,... T/\tau\in\N$.

It might be easily observed from \eqref{time-space-disc}, that in terms of the auxiliary variable $v_{\tau}^{k}$ which gives the equilibrium stress, the problem has the standard form of a linear elastic one and therefore could be numerically solved using any standard numerical procedure. However, BEM seems to be a natural choice, especially if we consider the case of zero body forces $f{=}0$, which we adopt for the rest of this work.

What is actually computed by BEM is the auxiliary field $v_\tau^k$, while we update  the elastic field $u_\tau^k$ by  \eqref{eqphysicalfield}, keeping in mind that $u_\tau^{k-1}$ is already known value at time step $k$. It is also important to notice that transformation   \eqref{eqauxfield}  appears also in the boundary condition on $\Gdir$, see \eqref{time-space-disc_02}, while tractions on $\GNeu$  are equal to tractions in the original visco-elastic problem, as shown in \eqref{time-space-disc_03}.

Taking into account the above explanation, the Somigliana displacement identity for the auxiliary variable $v^k$ can be written as
\begin{align}
\hspace{-2em} C(\xi)v^k_\tau(\xi){+}\int_{\Gamma} \!\!\!\!\!\!
-\ v^k_{\tau}(x)T(x,\xi)\,\d S_x{=}\int_{\Gamma} \mathfrak{t}(e(v^k_{\tau}))(x)U(x,\xi)\,\d S_x,
\end{align}
where, the weakly and strongly singular integral kernels $U(x,\xi)$ and $T(x,\xi)$ are the usual Kelvin fundamental solutions in displacements and tractions (two-point tensor fields) \cite{PaCa1997}, $C(\xi)$ is the coefficient tensor of the free term \cite{Mantic1993CM}, and the first integral represents the Cauchy principal value.

\subsection{Extension to multi-domain problems}\label{MultiDomail}
In problems of several bodies, where some of them may be visco-elastic or
merely elastic, we need to consider compatibility of displacements and
tractions equilibrium at common interfaces. Special attention is needed since,
while we solve the BEM system with respect to the auxiliary field $v_\tau^k$,
compatibility of displacement has to be considered for the displacement field
$u_\tau^k$. Thus, at the interface between two visco-elastic solids  $\Omega_1$
and $\Omega_2$ with relaxation times $\chi_1$ and  $\chi_2$, respectively,
the compatibility of displacements writes as
\begin{align}\label{dispComp}
&u_\tau^{k,1}=u_\tau^{k,2}\ \ \Rightarrow\ \
\frac{\tau}{\tau{+}\chi_1}v_\tau^{k,1}+\frac{\chi_1}{\tau{+}\chi_1}u_\tau^{k-1,1}
=\frac{\tau}{\tau{+}\chi_2}v_\tau^{k,2}+\frac{\chi_2}
{\tau{+}\chi_2}u_\tau^{k-1,2},
\end{align}
where a variable $q^{k,i}_\tau$ refers to the domain $\Omega_i$ at the $k^{\rm th}$
time step. For the case of elastic solids, where $\chi{=}0$,
eq.~\eqref{dispComp} cast to the usual equation considered in a BEM
formulation, that is $u_\tau^{k,1}{=}u_\tau^{k,2}$ reduces to
$v_\tau^{k,1}{=}v_\tau^{k,2}$, as in this case the auxiliary field $v_\tau^{k}$
obviously coincides with the displacement field $u_\tau^{k}$. Equilibrium of
tractions is considered for the total stress field defined in~\eqref{Kelvin} and
consequently for the tractions $\mathfrak{t}$ that correspond to the auxiliary
field $v_\tau^{k}$ and these tractions are directly computed in the BEM
formulation,
\begin{align}\label{tractionEqulibrium}
\mathfrak{t}^1(e(v_\tau^k))=-\mathfrak{t}^2(e(v_\tau^k)).
\end{align}
\subsection{Extension to contact problems utilizing the energetic approach in BEM}\label{ViscousContact}
Visco-elastic frictionless contact problems are numerically handled usually by utilizing FEM, cf.\ \cite{FEMVC1, FEMVC2, FEMVC3, FEMVC4, FEMVC5, FEMVC6}.
 To our best knowledge,
except for the specific case of rolling contact \cite{rolling1}, it is
the first time that a BEM formulation for contact problems of visco-elastic
solids is presented and fully explored, although it has been also used in \cite{RoPaMa_Visc} and originally proposed in \cite{Ro_AdhVisc}.
In order to solve the unilateral and/or adhesive contact problem of an assemblage of solids under (possible) contact to each other and/or  some outer rigid obstacles, we follow the general framework of energetic approaches to contact problems using BEM,
as it  is introduced in \cite{ECBEM}. \COL{Under this framework, the
minimization of the   potential  energy,    defined here in terms of the auxiliary variable $v_\tau^k$ from \eqref{eqauxfield},
\begin{align}\label{Gibbs-engr}
\calG(k\tau,v_{\tau}^{k})=\int_{\Omega} \frac12\bbC e(v_\tau^k){:}e(v_\tau^k) \dd x
-\int_{\GNeu}\!\!g_\tau^k{\cdot}v_\tau^k\dd S,
\end{align}
is required.}
The same procedure has also been utilized in \cite{RoPaMa_Visc}, however
without a  detailed presentation and  numerical testing of the BEM
formulation for   visco-elastic problems.

Here we assume a non-empty $\GC$ and write  the discretized condition \eqref{timedisc_04} in the form
\begin{align}
\label{time-space-disc_04}
&&&v_\tau^k{\cdot}\vec{n}\le -\frac\chi\tau u_\tau^{k-1}{\cdot}\vec{n},\ \ \ \
\mathfrak{t}_{\rm n}(e(v_\tau^k))\le0,\ \ \ \
(v_\tau^k{\cdot}\vec{n})\mathfrak{t}_{\rm n}(e(v_\tau^k))=0,\ \ \ \ \mathfrak{t}_{\rm t}(e(v_\tau^k))=0
&&\text{on }\GC,
\end{align}
which completes the system of equations \eqref{time-space-disc}.
Following the energetic approach in BEM, we obtain a convex
minimization  problem in terms of the auxiliary field $v_\tau^k$, in particular we have to solve the
quadratic-programming problem:
\begin{align}\label{LQ-programme}
&\left.\begin{array}{ll}
\text{minimize}&\displaystyle{
\calG(k\tau,v_{\tau}^k)
}
\\[.3em]
\text{subject to}&
v_\tau^k {\cdot}\vec{n}\le -
\text{\large$\frac\chi\tau$}
u_\tau^{k-1} {\cdot}\vec{n} \quad \text{on }\GC
\end{array}\right\}
\end{align}
with $\calG$ from \eqref{Gibbs-engr}. It is important to realize here that
in the quadratic-programming problem only the part of the auxiliary field
defined on $\GC$ represent active variables in the minimization procedure
\cite{ECBEM,PaMaRo13BEMI}.

Since the auxiliary variable $v_{\tau}^{k}$ gives the equilibrium stress, in contrast to the elastic field $u_{\tau}^{k}$, the domain integral appeared in $\calG$, under the assumption of zero body forces, can be expressed as a boundary one  through the so-called Clapeyron theorem, i.e.\
\begin{align}
\int_{\Omega}\!\frac12\bbC e(v_{\tau}^{k}){:}e(v_{\tau}^{k})\dd x=\frac12\int_{\Gamma}\!\mathfrak{t}(e(v^k_\tau))\cdot v_{\tau}^{k}\dd x,
\end{align}
and finally the stored energy in terms of $v_{\tau}^{k}$ and in a boundary form, that we
have to minimize, is given as
\REM{Tom: I used the boundary integral here - is it OK?}
\begin{align}
\calG(k\tau,v_{\tau}^k)=
\frac12\int_{\Gamma}\!\mathfrak{t}(e(v^k_\tau))\cdot v^k_\tau\dd x
-\int_{\GNeu}\!\!g^k_\tau{\cdot}v^k_\tau\dd S
\label{TotalPotentialEnergyAuxBoundary}
\end{align}
for which, standard techniques presented in \cite{ECBEM} might be used  to
numerically handle the above minimization problem by utilizing BEM.

With the above strategy, we estimate the sum of the stored elastic
energy  and the dissipated energy. However, we sometimes are interested in
visualizing the spatial distribution of the accumulated dissipated energy
due to viscosity, that is the term
$\int_0^t\!\chi\bbC e(\DT u){:}e(\DT u)\,\d t$ in~\eqref{TotalEngr}.
This is
meaningful for the vast majority of visco-elastic problems, and not only for
  contact problems we study in this section. It is a standard procedure in BEM,
that after solving the boundary value problem we  compute
displacements as well as   stresses and strains in the whole
domain by using the boundary values of displacements and tractions \cite{PaCa1997}. Having computed
 the stress and strain tensors in the required internal points
 for any time  $t_k$, we may easily compute the above time integral for any time  by using  the previous time history.
\section{Other linear visco-elastic rheologies}\label{sect-other}
The above method can be modified for other rheologies assuming again
like in \eqref{Kelvin} that all the viscous and the elastic responses have
the same tensorial character and thus are fully described just by only one
tensor and several scalar constants. A generalized linear visco-elastic model, consisting of an assemblage of the   Maxwell and
Kelvin-Voigt  elements together with free springs and dampers in series and/or parallel, might be represented by the following constitutive stress-strain relation in the  form of a differential equation  \cite{Brinson}:
\begin{align}\label{eq:GDE}
\sum_{k=0}^n \xi_k\frac{d^k\sigma}{dt^k}=\bbC e\left(\sum_{k=0}^m\chi_k\frac{d^k u}{dt^k}\right).
\end{align}

Obviously, certain restrictions on coefficients $\chi_k$ and $\xi_k$ exist,
see a detailed discussion in \cite{Flugge}.

Let us briefly present only a few special
cases for which all the manipulation can lucidly be demonstrated and which
simultaneously cover rheological models standardly used in most applications.
Nevertheless, we could routinely continue for more complex rheologies with
higher-order time derivatives on both sides, but the algebraic
manipulation would become complicated and the requirement for an
equal-tensorial character more restrictive. For simplicity,
 in this section  we do not consider the unilateral contact,
i.e.\ $\GC=\emptyset$, and, like before, we neglect inertial and external
bulk forces.
%
We  further   restrict ourselves, for implementation   and notational
purposes, to the case of the second-order stress-strain relation in \eqref{eq:GDE}, which for $n=m=2$
 is  given in the following form:
\begin{align}\label{eq:GDE2}
\xi_2\DDT\sigma+\xi_1\DT\sigma+\xi_0\sigma=
\bbC e(\chi_2 \DDT u+\chi_1 \DT u+\chi_0 u),
\end{align}
requiring some initial conditions for displacements and stresses and
their time derivatives of at most of the first order, depending on
the values of parameters $\chi_k$ and $\xi_k$.
The general form of equations that governs the system is,
\begin{subequations}\label{difsystem}
\begin{align}\label{difsystem-a}
&&&\mathrm{div}\,\sigma=0 \ &&\text{on }\Omega,&&&&\\
&&&u=w &&\text{on }\GDir,\\
&&&\sigma\vec{n}=g &&\text{on }\GNeu.
\end{align}
\end{subequations}

The implicit  time discretisation of eq. \eqref{eq:GDE2} assuming a fixed time step $\tau$, leads to
\begin{align}
\xi_2\frac{\sigma_\tau^k{-}2\sigma_\tau^{k-1}{+}\sigma_\tau^{k-2}}{\tau^2}+
\xi_1\frac{\sigma_\tau^k{-}\sigma_\tau^{k-1}}\tau+\xi_0\sigma_\tau^k
=\bbC e\bigg(\chi_2\frac{u_\tau^k{-}2u_\tau^{k-1}{+}u_\tau^{k-2}}{\tau^2}
+\chi_1\frac{u_\tau^k{-}u_\tau^{k-1}}\tau+\chi_0 u_\tau^k\bigg)
\end{align}
and, after an elementary algebra, the time-discrete variant of
\eqref{eq:GDE2} and \eqref{difsystem} reads as
\begin{align}
\mathrm{div}\,\sigma_\tau^k &=0\ \ \ \text{ with }\nonumber \\
\sigma_\tau^k
&=\bbC e\bigg(\frac{\chi_2{+}\tau\chi_1{+}\chi_0\tau^2}
{\xi_2{+}\tau\xi_1{+}\xi_0\tau^2}u_\tau^k
-\frac{2\chi_2{+}\tau\chi_1}{\xi_2{+}\tau\xi_1{+}\xi_0\tau^2}u_\tau^{k-1}
+\frac{\chi_2}{\xi_2{+}\tau\xi_1{+}\xi_0\tau^2}u_\tau^{k-2}\bigg)\nonumber
\\&\qquad
+\frac{2\xi_2{+}\tau\xi_1}{\xi_2{+}\tau\xi_1{+}\xi_0\tau^2}\sigma_\tau^{k-1}-
\frac{\xi_2}{\xi_2{+}\tau\xi_1{+}\xi_0\tau^2}\sigma_\tau^{k-2} \qquad \text{on} \quad \Omega,
\end{align}
 completed by the boundary conditions  $u_\tau^k=w_\tau^k$  on $\GDir$ and $\sigma_\tau^k\vec{n}=g_\tau^k$  on  $\GNeu$.

 The
implementation of BEM relies on $\mathrm{div}\,\sigma_\tau^{k-1}=0$ and
$\mathrm{div}\,\sigma_\tau^{k-2}=0$, and furthermore, likewise in
\eqref{eqauxfield}, on the definition of an auxiliary field of the general
form
\begin{align}\label{eqauxfield-g}
v_\tau^k=\frac{\chi_2{+}\tau\chi_1{+}\chi_0\tau^2}
{\xi_2{+}\tau\xi_1{+}\xi_0\tau^2}u_\tau^k
-\frac{2\chi_2{+}\tau\chi_1}{\xi_2{+}\tau\xi_1{+}\xi_0\tau^2}u_\tau^{k-1}
+\frac{\chi_2}{\xi_2{+}\tau\xi_1{+}\xi_0\tau^2}u_\tau^{k-2},
\end{align}
giving
\begin{equation}\label{sigmaktau}
\sigma_\tau^k
=\bbC e(v_\tau^k)+\frac{2\xi_2{+}\tau\xi_1}{\xi_2{+}\tau\xi_1{+}\xi_0\tau^2}\sigma_\tau^{k-1}-
\frac{\xi_2}{\xi_2{+}\tau\xi_1{+}\xi_0\tau^2}\sigma_\tau^{k-2}
\qquad \text{on} \quad \Omega.
\end{equation}

The transformed system of equations that we actually solve using BEM has the
form
\begin{subequations}\label{BEMsystem}
\begin{align}\label{BEMsystem_1}
&&&\mathrm{div}\,\bbC e(v_\tau^k)=0\ \ \ \
\quad \text{on }\Omega,\\
&&&v_\tau^k =\frac{\chi_2{+}\tau\chi_1{+}\chi_0\tau^2}{\xi_2{+}\tau\xi_1{+}\xi_0\tau^2}w_\tau^k
-\frac{2\chi_2{+}\tau\chi_1}{\xi_2{+}\tau\xi_1{+}\xi_0\tau^2}w_\tau^{k-1}
+\frac{\chi_2}{\xi_2{+}\tau\xi_1{+}\xi_0\tau^2}w_\tau^{k-2} \quad \text{on }\GDir,\\
&&&\mathfrak{t}(e(v_\tau^k))=
\big(\bbC e(v_\tau^k)\big)\big|_{\Gamma}\vec{n}  =g_\tau^k-\frac{2\xi_2{+}\xi_1\tau}{\xi_2{+}\xi_1\tau{+}\xi_0\tau^2}g_\tau^{k-1}+\frac{\xi_2}{\xi_2{+}\xi_1\tau{+}\xi_0\tau^2}g_\tau^{k-2} \quad \text{on }\GNeu.
\end{align}\end{subequations}
Solving the above system with BEM we obtain the pair $v_{\tau}^k$ and $\mathfrak{t}(e(v_\tau^k))$, for each time step $k$. Then, we may also compute $\sigma_{\tau}^k$, by evaluating $\bbC e(v_\tau^k)$ in $\Omega$ by standard BIR \cite{PaCa1997} and   adding $\sigma^{k-1}$ and $\sigma^{k-2}$ according to  \eqref{sigmaktau}. The reconstruction of the physical displacement field is carried out by solving eq.~\eqref{eqauxfield-g} for $u_\tau^k$,
\begin{equation}\label{utauk}
u_\tau^k=\frac{\xi_2{+}\tau\xi_1{+}\xi_0\tau^2}{\chi_2{+}\tau\chi_1{+}\chi_0\tau^2}v_\tau^k
+\frac{2\chi_2{+}\tau\chi_1}{\chi_2{+}\tau\chi_1{+}\chi_0\tau^2}u_\tau^{k-1}
-\frac{\chi_2}{\chi_2{+}\tau\chi_1{+}\chi_0\tau^2}u_\tau^{k-2},
\end{equation}
and, following \eqref{sigmaktau}, the total traction (physical traction) field  $p_\tau^k=\sigma_{\tau}^k\vec{n}$ is  reconstructed by
\begin{align}\label{tracRecon}
p_\tau^k=\mathfrak{t}(e(v_\tau^k)) +\frac{2\xi_2{+}\xi_1\tau}{\xi_2{+}\xi_1\tau{+}\xi_0\tau^2}p_\tau^{k-1}-\frac{\xi_2}{\xi_2{+}\xi_1\tau{+}\xi_0\tau^2}p_\tau^{k-2}.
\end{align}

All the necessary initial values, appearing above for discrete time lower than zero,
are assumed to be equal to zero. Calculation of characteristic physical parameters
is just a post-processing procedure and depends on each specific model. E.g., elastic
stresses of the Kelvin-Voigt model can be obtained recursively by applying the
elastic stress operator $\bbC e(\cdot)$ to \eqref{eqphysicalfield}, which is a
particularization of \eqref{utauk}. Some of the models that could be represented by the
second order differential equation \eqref{eq:GDE2} are listed in Table \ref{tab:models};
see also \cite{Brinson}. It is worth mentioning that the system of equations
\eqref{BEMsystem} could obviously be solved by any other appropriate numerical
method (e.g., FEM as in \cite{MesquitaCoda2001}), and that more complicated visco-elastic
models of a higher-order, i.e.\ $m>2$ or $n>2$ in \eqref{eq:GDE}, could be accomplished
within the current framework with the only difference that higher order derivatives will appear.

Within the class of constitutive relations defined by \eqref{eq:GDE2}, we will
consider several selected rheologies shown in Table \ref{tab:models}.  The
discrete energy estimates like \eqref {TotalEngr+disc}
can be derived for each of them
after suitable, sometimes rather complicated manipulation (not performed
in this article, however).

\def\NEWcheckmark{^{^{\text{\normalsize$\checkmark$}}}}
\def\NEWtimes{^{^{\text{\normalsize$\times$}}}}

\begin{table}
\caption{Some models that could be represented by the constitutive differential
equation \eqref{eq:GDE2} with pertinent coefficients $\chi$ and $\xi$,
present ($\neq 0$) indicated by $\checkmark$ or absent ($=0$) by $\times$.}
\medskip
\centering
    \begin{tabular}{ | c | p{13em} | l | l | l | l | l | l |}
 \hline
    Model $^{^{^{^{}}}}_{_{_{_{}}}}$ & Name & $\chi_0$ & $\chi_1$ & $\chi_2$ & $\xi_0$ & $\xi_1$ & $\xi_2$
   \\ \hline\hline
    \includegraphics[scale=0.25]{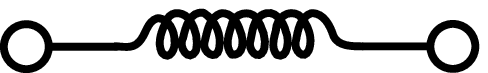} $^{^{^{^{}}}}_{_{_{_{}}}}$ &
      Elastic (Hooke) solid & $\checkmark$ & $\times$ & $\times$ & $\checkmark$ & $\times$ & $\times$
 \\ \hline
    \includegraphics[scale=0.25]{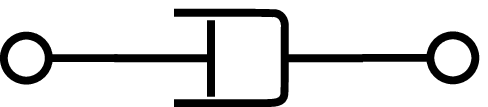} $^{^{^{^{}}}}_{_{_{_{}}}}$
     & Viscous (Newton) fluid & $\times$ & $\checkmark$ & $\times$ & $\checkmark$ & $\times$ & $\times$
  \\ \hline
    \includegraphics[scale=0.25]{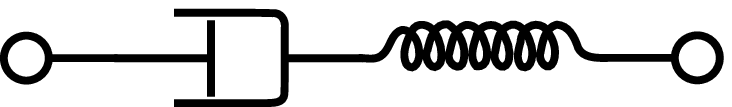} $^{^{^{^{}}}}_{_{_{_{}}}}$
     & Maxwell fluid & $\times$ & $\checkmark$ & $\times$ & $\checkmark$ & $\checkmark$ & $\times$
  \\ \hline
$\text{\includegraphics[scale=0.25]{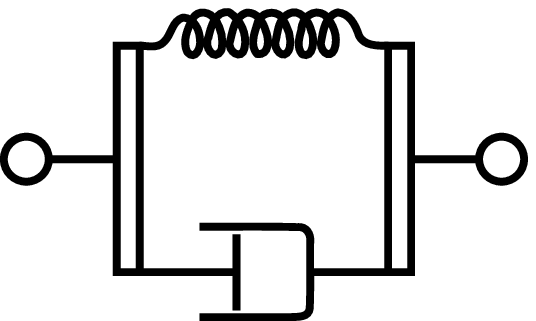}}^{^{^{^{}}}}$ &
     $^{^{\text{\normalsize Kelvin-Voigt solid}}}$ & $\NEWcheckmark$
    & $\NEWcheckmark$ & $\NEWtimes$ & $\NEWcheckmark$ & $\NEWtimes$ & $\NEWtimes$
  \\ \hline $\text{\includegraphics[scale=0.25]{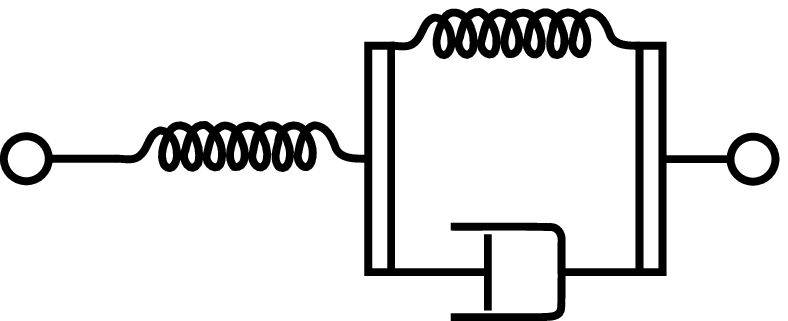}}^{^{^{^{}}}}$
   & $\!\!^{_{\begin{array}{l}\text{\normalsize Boltzmann or Standard linear}\\[-.2em]\text{\normalsize or 3-parameter solid}\end{array}}^{^{}}}$
    & $\NEWcheckmark$ & $\NEWcheckmark$ & $\NEWtimes$ & $\NEWcheckmark$ & $\NEWcheckmark$ & $\NEWtimes$
  \\ \hline
$\text{\includegraphics[scale=0.25]{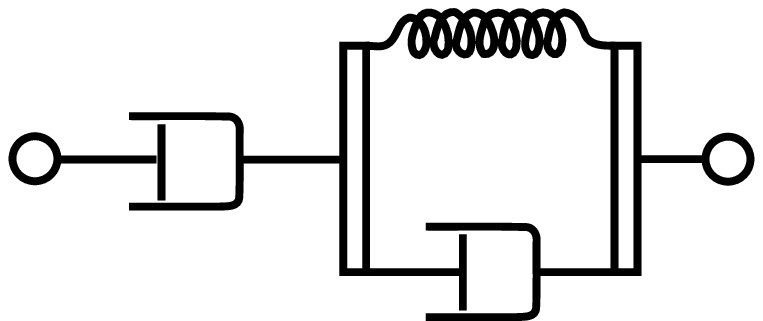}}^{^{^{^{}}}}$
    & $^{^{\text{\normalsize Jeffreys or 3-parameter fluid}}}$ & $\NEWtimes$ &
    $\NEWcheckmark$ & $\NEWcheckmark$ & $\NEWcheckmark$ & $\NEWcheckmark$ & $\NEWtimes$
  \\ \hline
$\text{\includegraphics[scale=0.25]{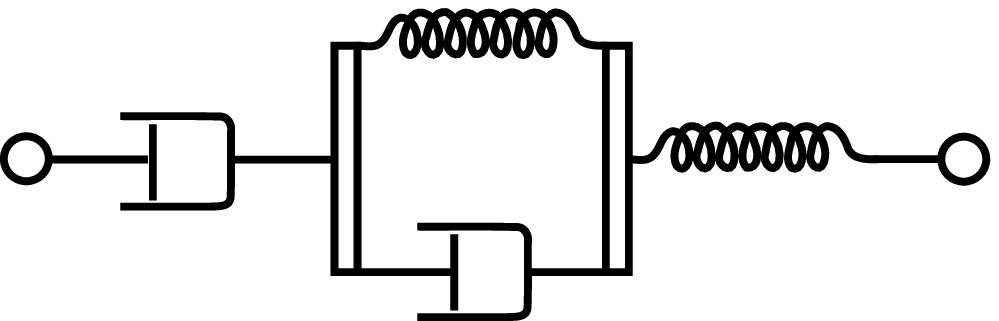}}^{^{^{^{}}}}$
    & $^{^{\text{\normalsize Burgers or 4-parameter fluid}}}$ & $\NEWtimes$ & $\NEWcheckmark$
    & $\NEWcheckmark$ & $\NEWcheckmark$ & $\NEWcheckmark$ & $\NEWcheckmark$
   \\ \hline
$\text{\includegraphics[scale=0.25]{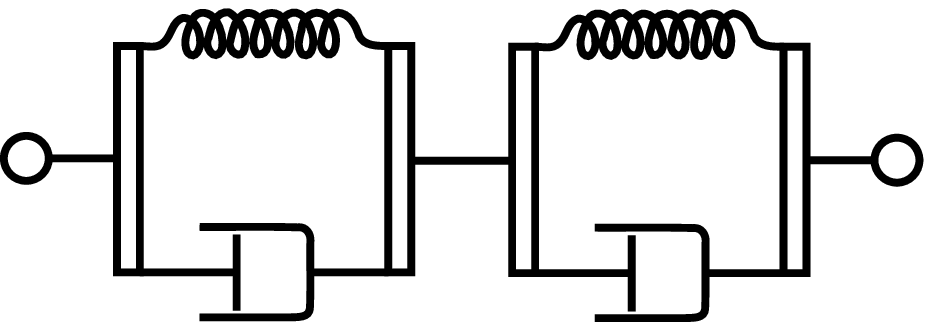}}^{^{^{^{}}}}$
    & $^{^{\text{\normalsize 4-parameter solid}}}$ & $\NEWcheckmark$ & $\NEWcheckmark$ &
    $\NEWcheckmark$ & $\NEWcheckmark$ & $\NEWcheckmark$ & $\NEWtimes$
   \\ \hline
    \end{tabular}
    \label{tab:models}
\end{table}
\section{Numerical examples}\label{Numericalexamples}
The above introduced framework has been implemented in an open BEM Java code \cite{OpenBem} with capabilities of 2D and 3D elastostatic analysis, among others. This code is supplied with all the necessary
``modules'' for the energetic approach in BEM used for contact problems,  and has also been employed in several related works of the authors~\cite{RoPaMa_Visc,ECBEM,PaMaRo13BEMI}.
\subsection{Visco-elastic creep behaviour}\label{Simple}\REM{Application00053}
This first example might be seen as a ``benchmark'', since it is one of the most frequent examples, met in the literature in order to compare numerical to analytical solutions of visco-elasticity (e.g. in \cite{MesCod02}).

\begin{figure}[ht!]
   \centering
   \def\svgwidth{.93\columnwidth}
   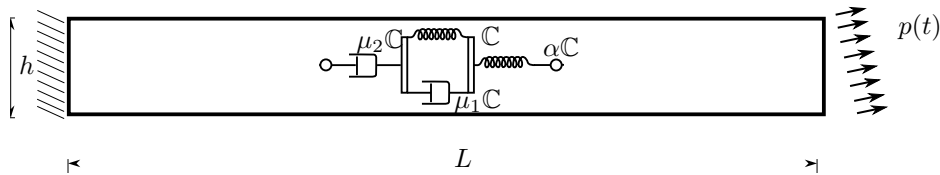
\vspace*{-1em}
   \caption{Geometry of the problem and physical interpretation.}
   \label{fig-geometry-prA}
\end{figure}

Figure~\ref{fig-geometry-prA} depicts the geometry and boundary conditions of the problem together with a physical interpretation of the visco-elastic mechanism of the material. The physical properties and the geometry of the problem are given in Table~\ref{table-geometry-and-properties}, where for the first variant of the problem we assume a Kelvin-Voigt material with viscosity $\mu_1$, without  the spring $\alpha\bbC$ and the damper $\mu_2\bbC$  depicted in Figure~\ref{fig-geometry-prA}. The uniform BEM mesh for this problem has 180 linear elements. Two time steps have been used, a coarse and a fine one, $\tau_c$=10\;(days) and $\tau_f$=1\;(day) respectively, in order to observe numerically the accuracy of the time integration scheme. Prescribed tractions on the right-hand side of the domain  have normal and tangential components $p_{\rm n}$=5\;(N/mm$^2)$ and $p_{\rm t}$=0, respectively. The total time of analysis is $T$=800\;(days). The external loading is removed at time $t_r$=400\;(days), i.e. after this time $p_n=0$.
\begin{table}[ht!]
 \caption{Elastic and geometrical properties of models used in Example A.}
 \medskip
 \centering
\begin{tabular}{l l}
 \hline
$L$ (mm)	&	800\\
$h$ (mm)	&	100\\
  \hline \\
$\mu_1$ (days) & 45.454545\\
$E$ (kN/mm$^2$)	&	11\\
$\nu$ 	&	0.0 \\
 \hline
\end{tabular}
\label{table-geometry-and-properties}
\end{table}
\begin{figure}[ht!]
\begin{center}
\includegraphics[width=.7\columnwidth]{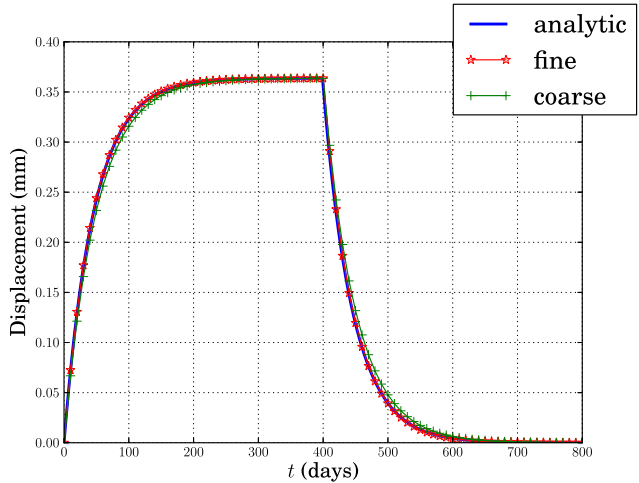}
\end{center}
\vspace*{-1.5em}
\caption{Displacement for the Kelvin-Voigt material, fine time partition solution shown here with one time point  per four steps.}
\label{fig-exA-disps}
\end{figure}

\begin{figure}[ht!]
\begin{center}
\includegraphics[width=.7\columnwidth]{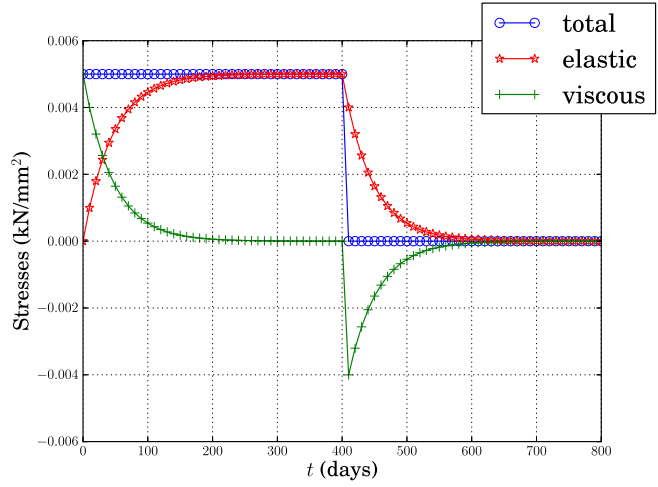}
\end{center}
\vspace*{-2em}
\caption{Stress $\sigma_{xx}$, for the Kelvin-Voigt material, at the centroid of the solid with one time point shown per ten steps, which means that only 80 time points are plotted, instead of   800 that actually have been calculated.}
\label{fig-exA-stresses}
\end{figure}

  Computed displacements are plotted in time in Figure~\ref{fig-exA-disps} together with the analytic solution, which can be  easily   deduced for this simple problem. Both numerical solutions for a coarse and a fine time step, are plotted. Notice that  the fine-time-step  solution is not shown in the plot for all time steps but only for those of the coarse partition of the time interval. An excellent  agreement of the fine-time-step  solution  with the analytic one is observed,  the coarse-time-step solution being also  very good. Figure~\ref{fig-exA-stresses} shows the evolution in time of the total stresses at the geometric center of the solid. Recall   that for the present case of the Kelvin-Voigt model, the \emph{total  stress field}, $\sigma_\tau^k =\bbC e(v_\tau^k)$,   corresponds directly to the auxiliary field $v_{\tau}^{k}$, while the \emph{elastic  stress  field}, $\bbC e(u_\tau^k)$, corresponds to the $u_{\tau}^{k}$ field. Then, the \emph{viscous stresses} can be computed as the difference of the total minus  elastic stresses.

In the second variant of this problem,   the prescribed tractions on the right-hand side   have components $p_{\rm n}$=0 and $p_{\rm t}$
=5\;(N/mm$^2)$, with the loading applied from the time
$t_i=80\;$(days) to $t_r=533.33$\;(days), while $T=800$\;(days). Numerical results are
obtained using time step $\tau$=1\;(day). In this case we show the spatial distribution of the dissipated energy density due to the
viscosity over the time interval $[0,T]$ for the Kelvin-Voigt model, and compare the kinematic  response of several visco-elastic rheologies presented in this article.

\begin{figure}[ht!]
\begin{center}
\includegraphics[width=.7\columnwidth]{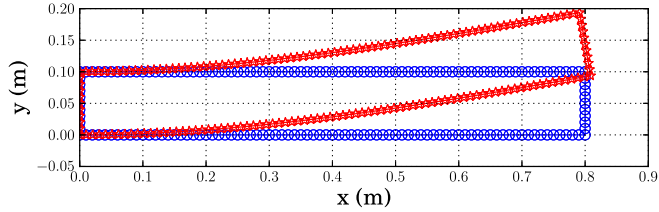}
\end{center}
\vspace*{-2em}
\caption{Deformed configuration, for the Kelvin-Voigt material, for the case of vertical loading, at time $t=T/2$.}
\label{fig-exA-disp2}
\end{figure}
Figure~\ref{fig-exA-disp2} shows the BEM mesh (used for both   variants of the problem) together with a deformed configuration for the case of vertical loading. In the next Figure~\ref{fig-exA-energy} the spatial distribution of the dissipated energy density $\int_0^T\!\chi\bbC e(\DT u){:}e(\DT u)\,\d t$, in (J/m$^2$), is visualized. It can be observed there, that the main part of the dissipated energy is accumulated, during the  evolution in time,  in a region close to the left fixed side of the solid where the highest normal stresses $\sigma_{xx}$ can be expected.
\begin{figure}[ht!]
\begin{center}
\includegraphics[width=.25\columnwidth, angle=-90]{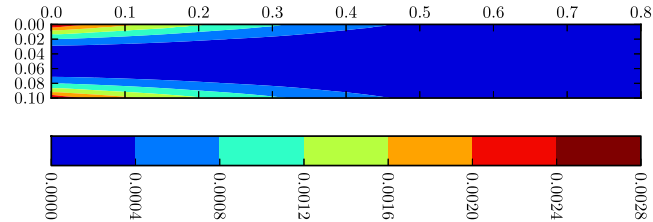}
\end{center}
\vspace*{-2em}
\caption{Spatial distribution of the dissipated energy density $\int_0^T\!\chi\bbC e(\DT u){:}e(\DT u)\,\d t$, in (J/m$^2$), for the case of vertical loading and the Kelvin-Voigt material.}
\label{fig-exA-energy}
\end{figure}

\begin{figure}[ht!]
\begin{center}
\subfigure[]{\includegraphics[width=.47\columnwidth]{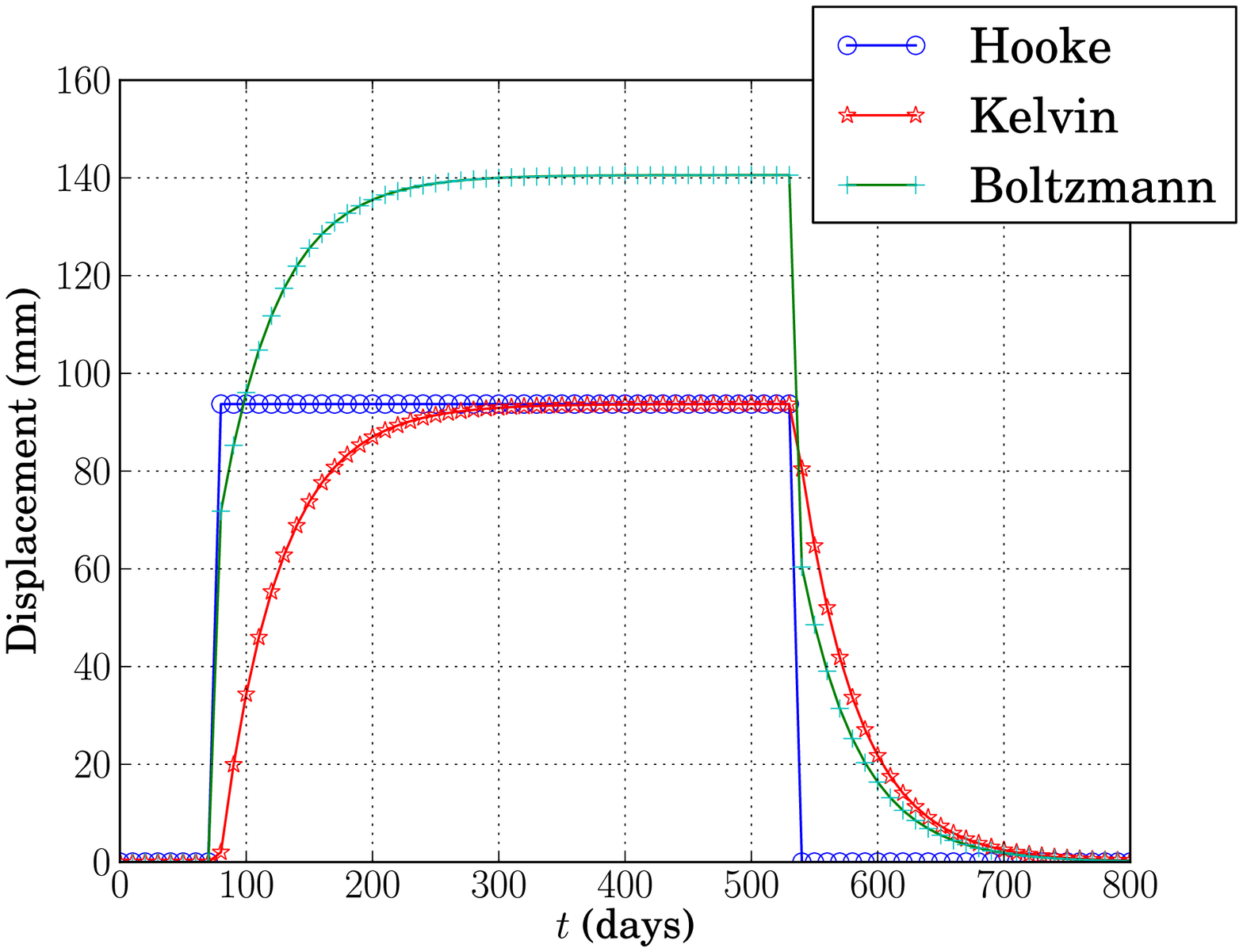}}
\ \ \
\subfigure[]{\includegraphics[width=.47\columnwidth]{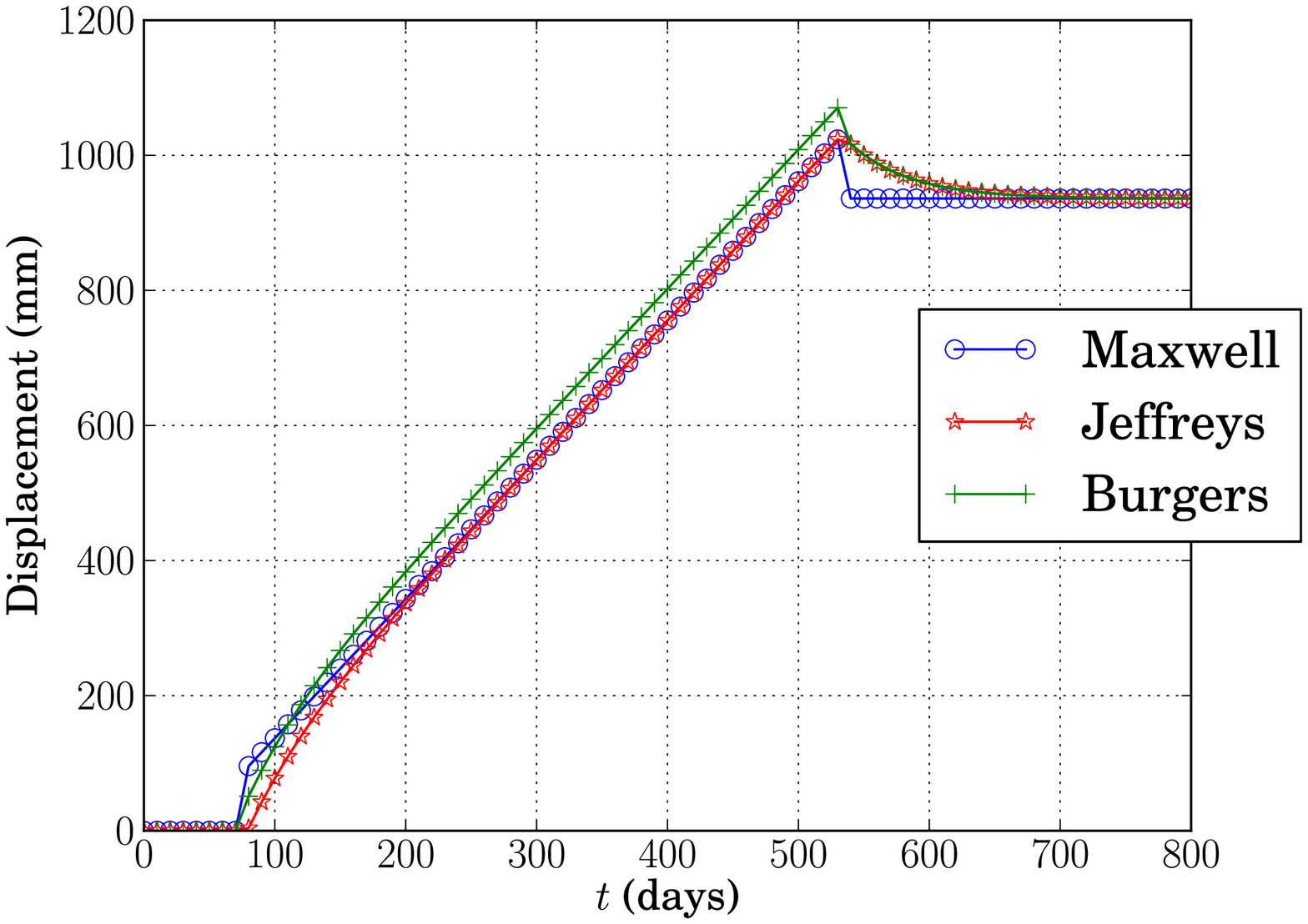}}
\end{center}
\vspace*{-2em}
\caption{Vertical displacement of the right-hand edge computed  by six different
rheology models, distinguished as (a) solid-type and (b) fluid-type.}
\label{fig:Compare}
\end{figure}

For the other   visco-elastic models studied we   use  the parameter values $\alpha{=}2$, $\mu_2{=}\mu_1$, where their nonzero values are required. For example, we may assume existence of the damper $\mu_2$ and the spring of stiffness $\bbC$, with simultaneous absence of the other two components, in order to simulate the Maxwell  model. The results are shown in Figure \ref{fig:Compare}, where models have been divided into two categories: (a) solid-type and (b) fluid-type, because of different order of response values. It might be observed in this figure the ability of the algorithm to compute a jump in displacement due to a jump of forces  for the case of both the Hooke and Boltzmann models  in contrast to the Kelvin-Voigt model, where a smoother increase of displacement takes place.
\subsection{3D analysis of an ellipsoidal cavity embedded in an infinite medium}\label{ThreeD}
\begin{figure}[ht!]
\begin{center}
\subfigure[]{\includegraphics[width=.4\textwidth]{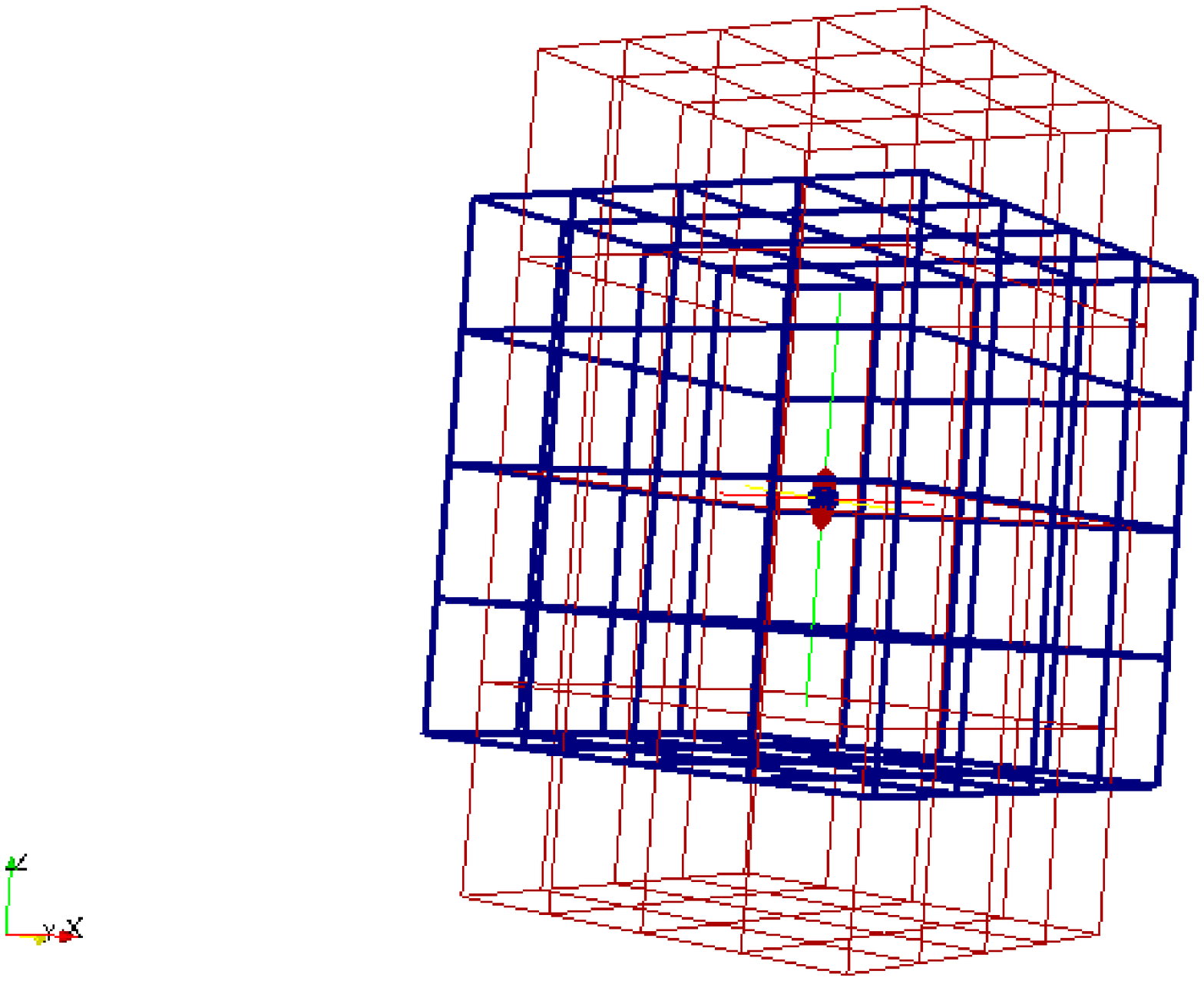}}
\subfigure[]{\includegraphics[width=.4\textwidth]{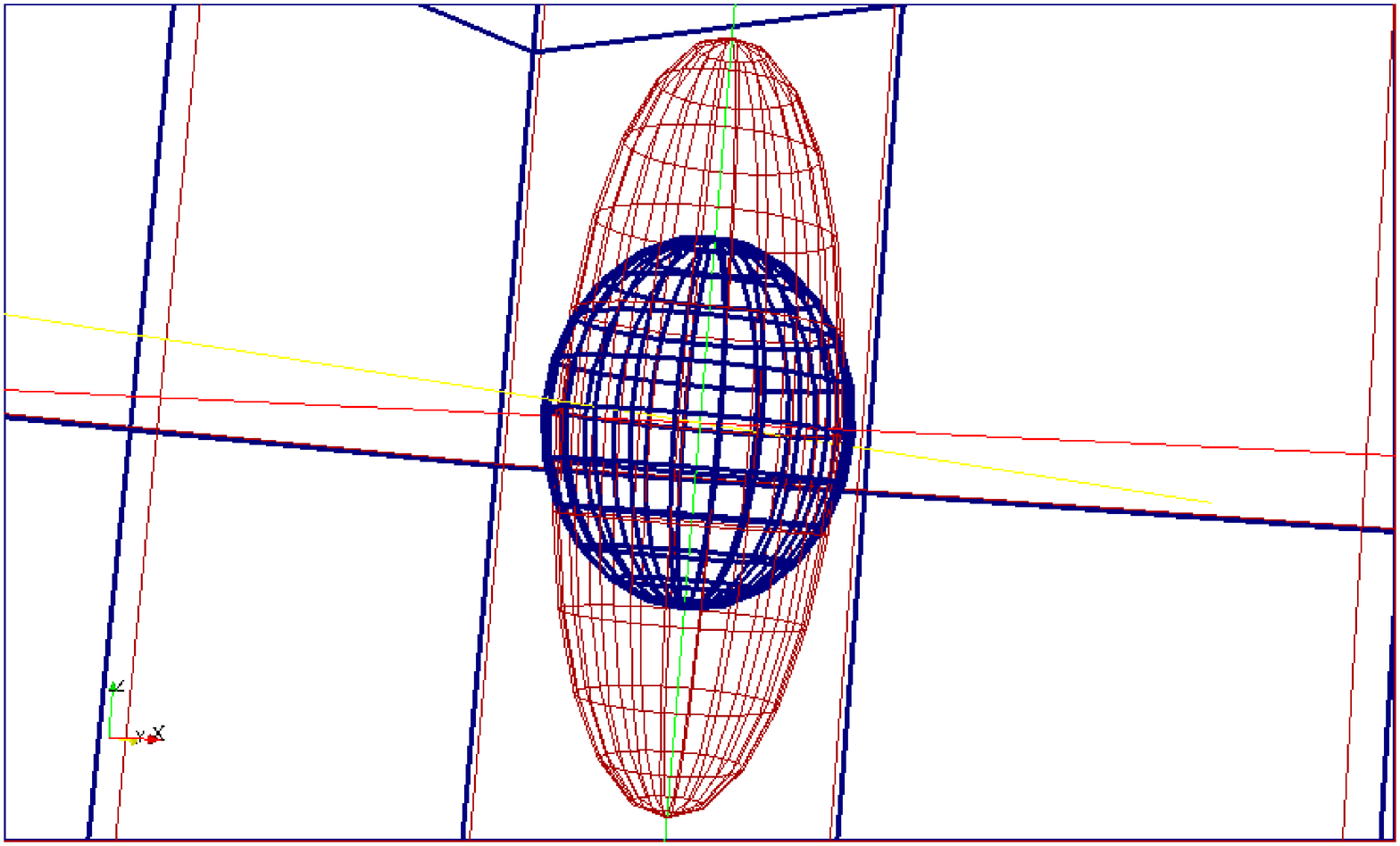}}
\end{center}
\vspace*{-2em}
\caption{(a) Undeformed and deformed BEM mesh of cube   with the embedded ellipsoidal  cavity, shown in detail in (b), at time $t=$400(sec). Scale factor of 500 is used   to magnify displacements.}
\label{fig:ex_2_1}
\end{figure}

\begin{figure}[ht!]
\begin{center}
\includegraphics[width=.7\columnwidth]{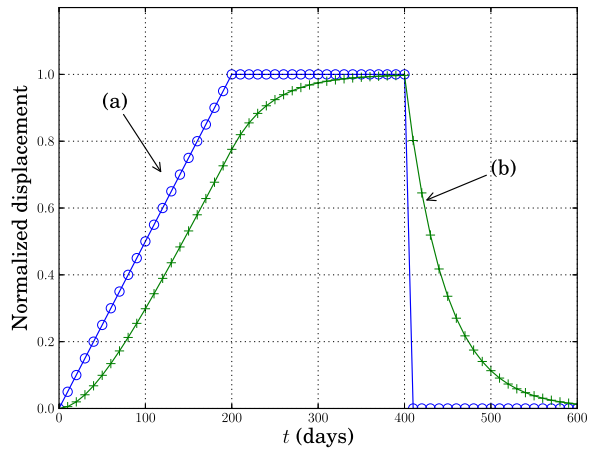}
\end{center}
\vspace*{-2.em}
\caption{Time evolution of the  displacement of the positive $Z$ pole of ellipsoid, normalized   by the maximum value of this displacement in the elastic case ($u_e^{max}=2.224$mm). (a) Elastic material, (b)  Visco-elastic Kelvin-Voigt material (the maximum value of this displacement is $u_v^{max}=2.218$mm). The time evolution pattern of the external loading   coincides with   the displacement evolution in the elastic  case.}
\label{fig-exB-results}
\end{figure}

\REM{Application00055}
 This example  shows the   capabilities of the procedure developed and implemented  also for   3D visco-elastic problems, see  \cite{SchanzAntes,Gao_Peng},  for other 3D BEM implementations. The  problem  of an ellipsoidal cavity in  a visco-elastic medium under remote stress field is solved. The Kelvin-Voigt material   considered  has Young's modulus $E$=70\;(GPa), Poisson's ratio $\nu$=0.35 and relaxation time $\chi$=45.454545\;(days). The remote stress field is applied on a cube with  side   length $L$=36\;(m) representing an infinite visco-elastic medium with an  embedded  ellipsoidal cavity placed in its center. The geometry of the ellipsoid is defined in Cartesian coordinates by the equation
\begin{align}
\frac{x^2}{a^2}+\frac{y^2}{b^2}+\frac{z^2}{c^2}=1,
\end{align}
with $a=$0.8\;(m), $b=$0.9\;(m) and $c$=1\;(m). The   BEM mesh of the ellipsoid  consists of 264 four node isoparametric quadrilateral elements, while the   cube boundary is discretised by 96 elements, see   Figure~\ref{fig:ex_2_1}.
Uniform normal tractions  $\sigma_x$=25\;(GPa), $\sigma_y$=25\;(GPa) and $\sigma_z$=100\;(GPa)  are applied    on the cube faces  perpendicular to the $x$-, $y$- and $z$-axis, respectively. The cavity boundary is free. The time pattern of  the load  has three parts:   initially the load increases linearly with time, then it is constant in time, and finally it jumps down to zero, as can be seen in Figure~\ref{fig-exB-results}. As only Neumann boundary conditions are prescribed, to avoid rigid body motions we apply the \emph{F1 method} of~\cite{BlasquezVlado}; to the best of our knowledge, first time implemented in the 3D case.
\subsection{Visco-elastic solid in contact}\label{Contact}
\REM{Application00054--half disc, Application00018-quarter disc}
A problem including frictionless contact between a viscoelastic solid and a rigid obstacle   is solved by the  BEM,   to the best of our knowledge, for the first time. The Kelvin-Voigt rheology is assumed.
\begin{figure}
   \centering
   \def\svgwidth{.9\columnwidth}
   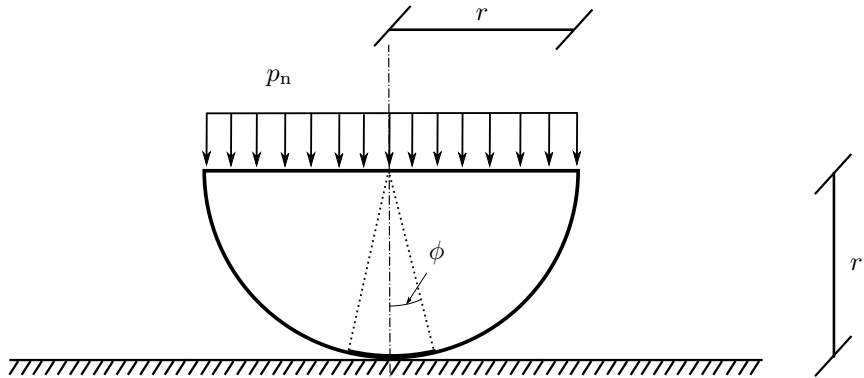
   \vspace*{-2em}
\caption{A visco-elastic half disk pressed against  the rigid foundation.}
\label{fig:geomB}
\end{figure}
In particular, the indentation of a half disk against a rigid  foundation   is considered under plane strain conditions. In this advancing contact problem the length of the contact zone depends on the load value. The problem  geometry is shown in Figure \ref{fig:geomB}. The radius of the disk is $r$=0.75m. The potential contact zone is defined by the angle $\phi$=13.5($^\circ$). Normal tractions are increased linearly in time from zero  to $p_{\rm n}$=-250\;(GPa) at time $t_p$=250\;(days) and then they are removed. We study the response up to the total time $T$=500\;(days). Tangential tractions along the whole straight edge of the half disk are zero. Due to the problem symmetry only the   quarter disc is modeled. The Kelvin-Voigt material   has Young's modulus $E$=70\;(GPa) and Poisson's ratio $\nu$=0.35.
For   comparison purposes,   three   relaxation times are considered:  $\chi=0$,
  $\chi=$45\;(days) and   $\chi=$22.5\;(days).
\begin{figure}[ht!]
\begin{center}
\subfigure[]{\includegraphics[width=.47\columnwidth]{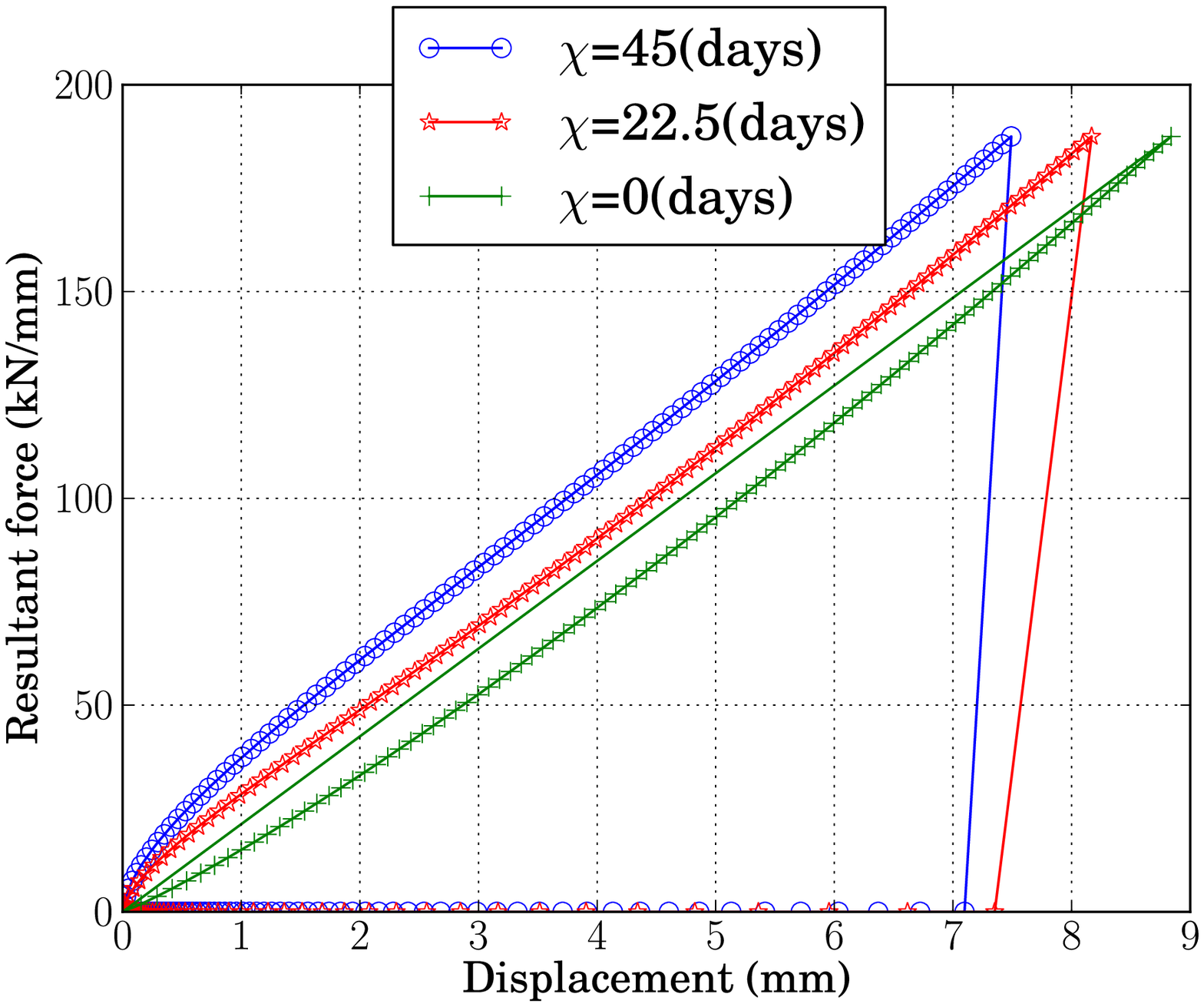}}
\ \ \
\subfigure[]{\includegraphics[width=.47\columnwidth]{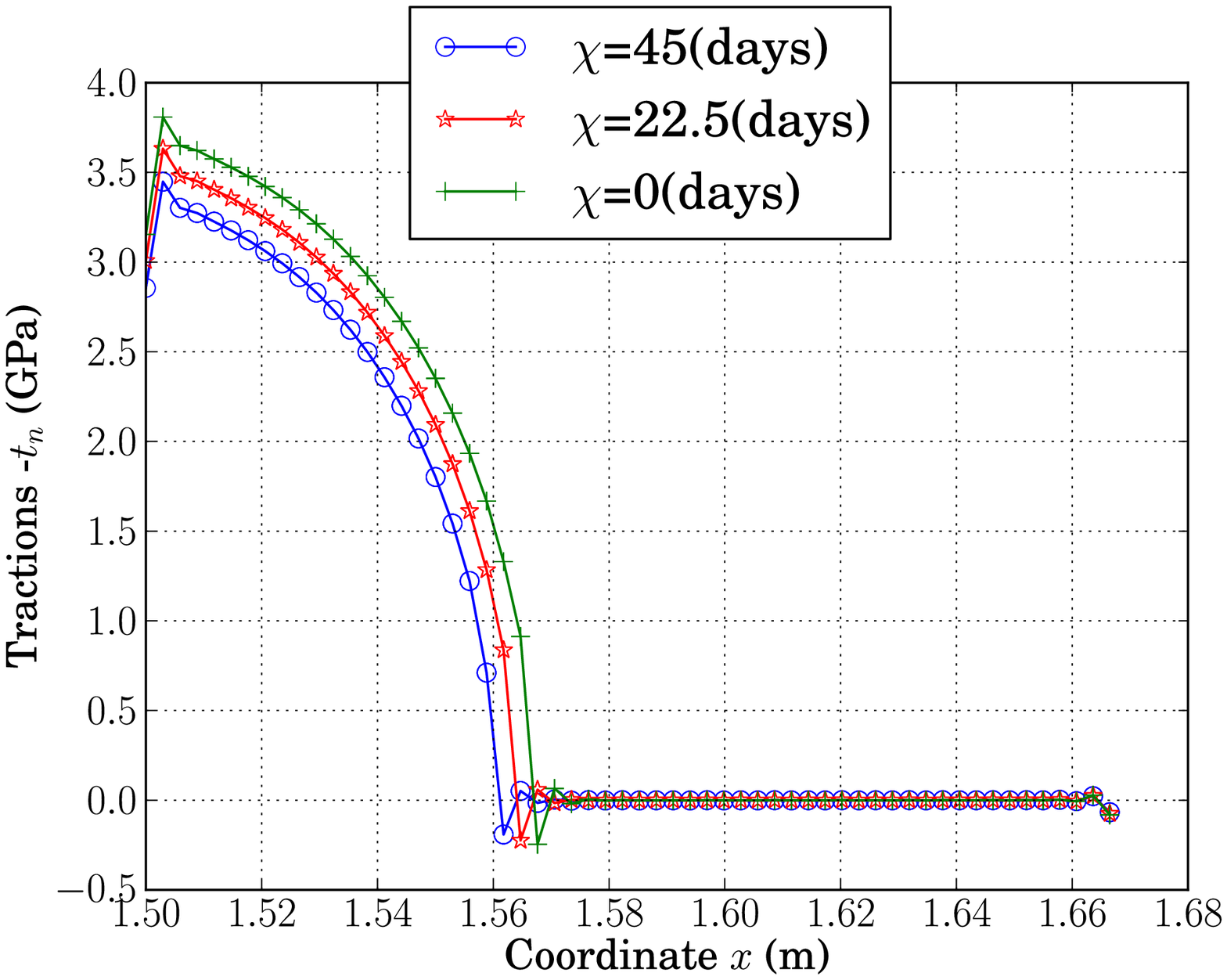}}
\end{center}
\vspace*{-2em}
\caption{(a) Total resultant vertical force on the horizontal side of the half disk versus the absolute value of the vertical displacement of the central point of this side.  (b) Normal elastic tractions along the possible contact zone at the time of peak loading $t_p$.}
\label{fig-exC_load_disp}
\end{figure}

 The numerical solution of this problem, which includes the determination of the contact zone, is accomplished as  described in Section~\ref{ViscousContact}, through the minimization of the potential energy.  The BEM mesh of the quarter disk consists of 270 linear elements with 60 elements along the possible contact zone defined by the angle $\phi$, 170 elements for the rest of the circular curve and 20 elements for each one of the two straight lines. The time step of $\tau$=2.5\;(days) is used, for the three relaxation times considered, resulting in 200 time steps.
\begin{figure}[ht!]
\begin{center}
\subfigure[]{\includegraphics[width=.47\columnwidth]{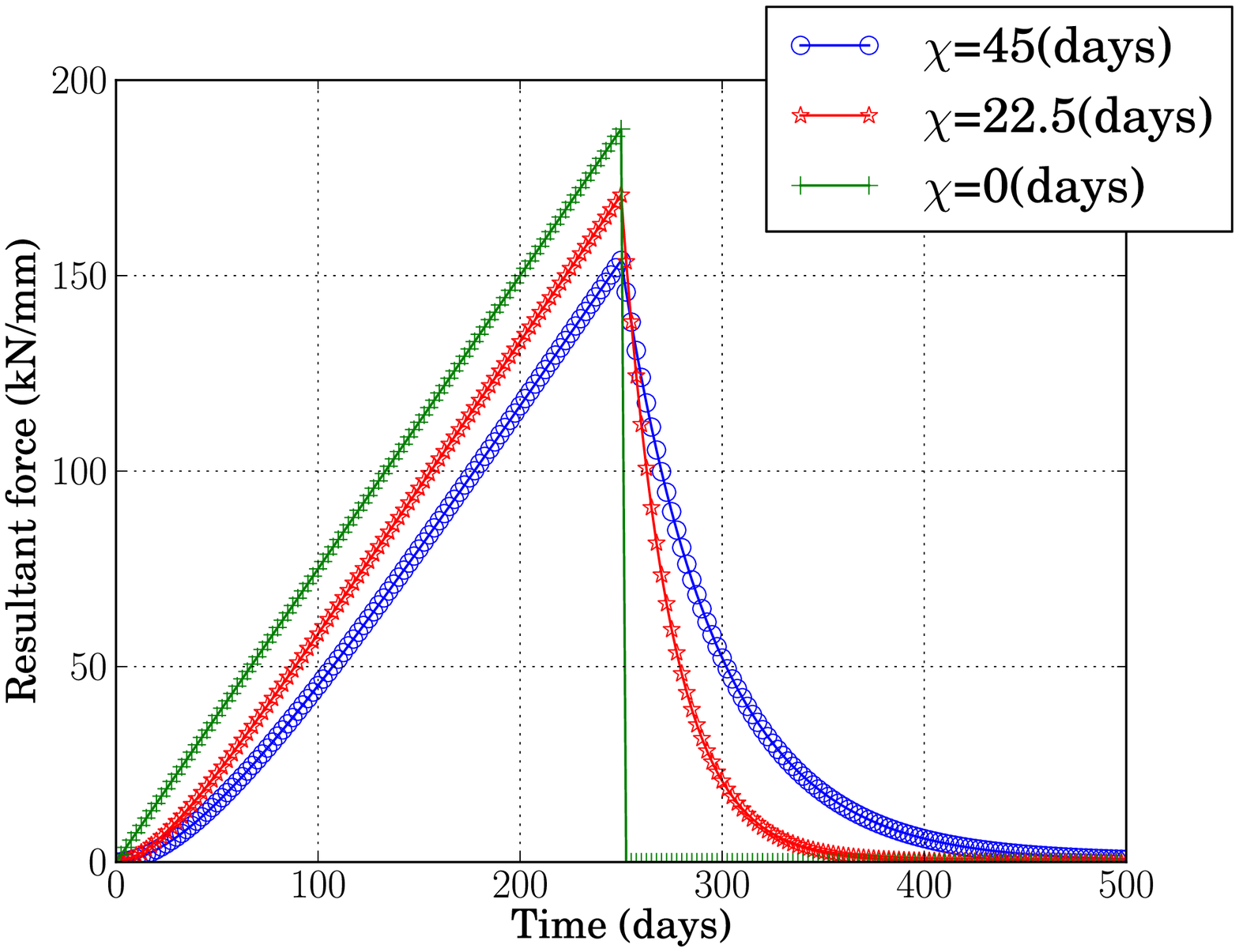}}
\subfigure[]{\includegraphics[width=.47\columnwidth]{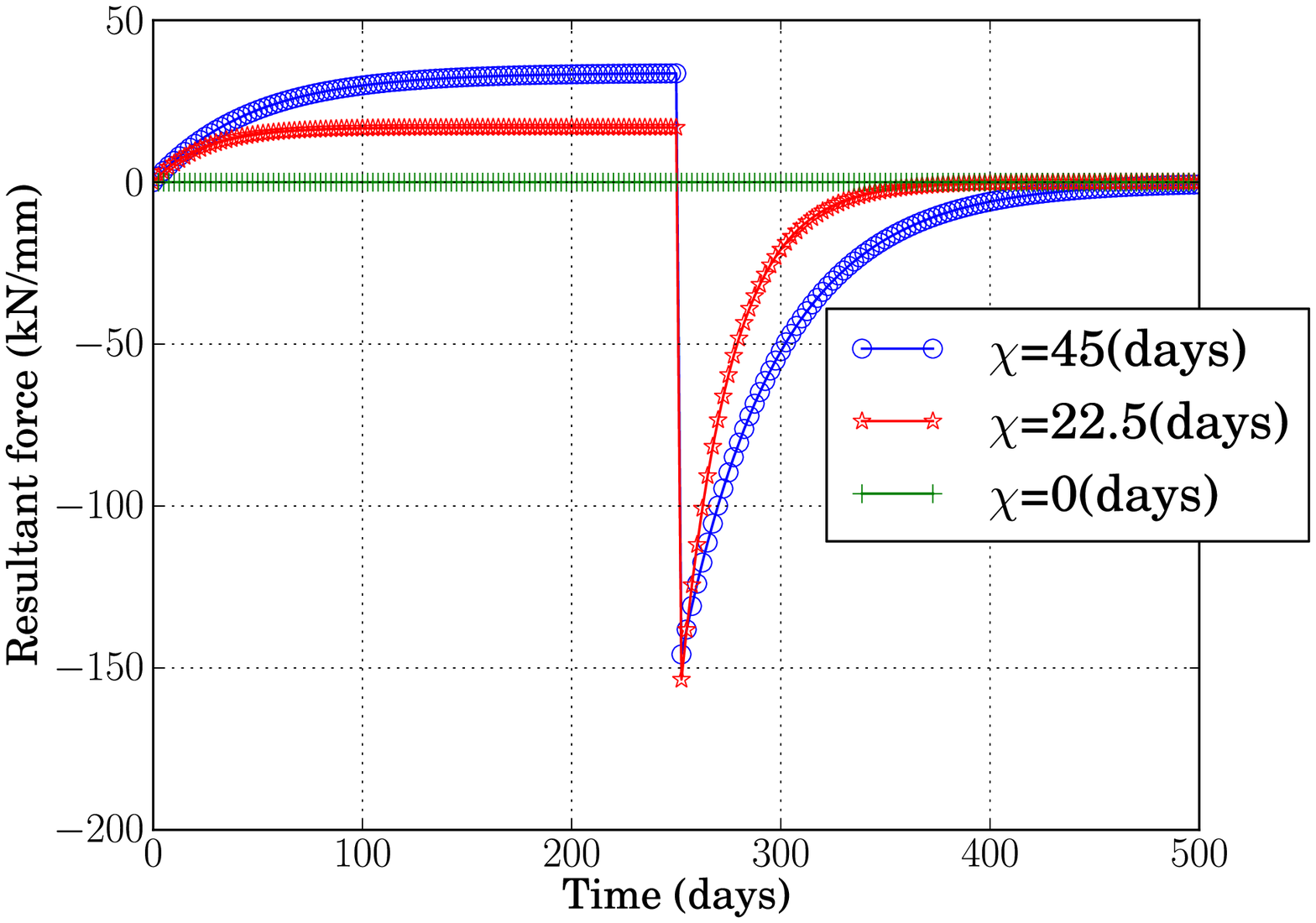}}
\end{center}
\vspace*{-2em}
\caption{(a) Elastic resultant force with time.
(b) Viscous resultant force with time.}
\label{fig-exC_load_time}
\end{figure}

The advancing contact problem is non-linear and this can be verified from Figure~\ref{fig-exC_load_disp}-(a), where for the non-viscous case and after the loading is removed the solution directly returns to the initial configuration. It might be seen there, that the straight line that connects the ``peak'' loading point back to the initial configuration is different  from the non-linear path   computed from the initial undeformed  configuration to the peak point. The behaviour of the visco-elastic cases is different, where we   observe that the greater the viscosity the greater the difference from the elastic case. For the viscous cases, we notice that after  the loading vanishes at the time $t_p$, the total force jumps to zero as well, while elastic and viscous forces of opposite signs still   remain and vanish progressively.  It is also easily verified  from Figure~\ref{fig-exC_load_disp}-(b)  that the length of the contact zone depends on $\chi$ value. It can been observed there that the greater the viscosity,   lower the length of the contact zone and   lower the maximum absolute value of the normal elastic tractions. This last observation may also be noticed in Figure~\ref{fig-exC_load_time}-(a), where the evolution in time of the elastic part of the resultant force is plotted for all three viscosity cases. Finally, in Figure~\ref{fig-exC_load_time}-(b)  the evolution of the viscous part of the resultant force is plotted, where it is interesting   to observe a jump and a finite peak in these viscous forces at the time of the loading removal $t_p$ for $\chi>0$.
\section{Conclusions}

\COL{In this paper, an advanced formulation for the solution of quasistatic
linear visco-elastic problems for a broad spectrum of rheologies, which
further develops the original proposal by Mesquita, Coda and co-workers
\cite{MesquitaCoda2001,MesCod02}, has been presented.}
%
The resulting problem can be solved using standard numerical methods such
as FEM and BEM.

We  have confined  ourselves to materials responding on the mechanical
loading in such a way that,
roughly speaking, the tensorial and the rheological features are separated;
this means only one tensor is used to describe all the elastic and viscous
processes which then are distinguished only be scalar constants.
Since, we have been  able to  cast the problem using boundary formulas
only,  and then  BEM appears as the most reasonable method in order to
solve both 2D and 3D problems.  After a certain
``computational cheap'' algebraic manipulation, only the standard
Kelvin's fundamental solution of elasticity
is required for the BEM
implementation. Furthermore, an extension and implementation to contact
problems of visco-elastic continua is presented  as well.

Using this formulation, the well known Kelvin-Voigt model has been scrutinized
and it has been shown that several other, more complex models, can be
confronted. A quite detailed presentation has been given for several models
using the Maxwell, Boltzmann, Jeffreys and Burgers rheologies.

Incorporation of this framework to existing BEM codes is very easy, at least for problems of visco-elasticity, since just a transformed auxiliary field has to be defined. After solving the problem for this auxiliary field, the actual stresses and displacements
can be easily reconstructed. For unilateral contact problems,
further features of the
energetic approach in BEM are needed. Numerical solutions of problems presented in this paper are accomplished by an in-house open BEM code, implemented in Java.

Some standard problems of 2D and 3D visco-elasticity as well as a problem of contact mechanics have been numerical solved and analysed in order to validate the suitability of the  methodology developed for solving realistic   visco-elasticity problems.

 An  extension of the current framework to problems of adhesive contact
or also to more complex problems, where interface damage and/or interface plasticity are taken into account,  is possible and into some extent has already been accomplished in other concurrent works of the authors (e.g. \cite{RoPaMa_Visc,KruPaRo_Visc}).
\section*{Acknowledgments}
{\small
The authors thank to two anonymous reviewers for their constructive comments, which were helpful
in improving the manuscript.
The authors acknowledge the support by the Junta de Andaluc\'{\i}a and Fondo Social Europeo
(Proyecto de Excelencia TEP-4051), by the Ministerio de Econom\'{\i}a y Competitividad (Proyecto MAT2012-37387),  as well as from the grants
201/09/0917, 201/10/0357, and 13-18652S  (GA \v CR) together with the institutional support RVO:\,61388998 (\v CR).}

\bibliographystyle{model3-num-names}
\bibliography{biblio}






\COL{
\appendix
\section*{Appendix: The energetics of selected rheological models}

All rheological models above allow for clear energetic balance, which is
important in many respects. We will illustrate it only for the standard
linear solid and, as special cases, for the Maxwell and the Kelvin-Voigt
models, i.e.\ \eqref{eq:GDE} for $m\le1$ and $n\le1$.

The energetics for the standard linear solid (and for Maxwell material too)
needs an introduction of one {\it internal variable} with the meaning
of a strain, let us denote it by $\pi$, acting in an additive
decomposition of the total strain $e(u)$, i.e.\
\begin{align}\label{e(u)=e-pi}
e(u)=e_\mathrm{el}+\pi.
\end{align}
The elastic strain $e_\mathrm{el}$ occurs on the ``serial'' elastic spring
(let us denote its elastic-moduli tensor by $\bbC_\mathrm{M}$)
while $\pi$ occurs on the ``parallel'' elastic spring
(with the elastic moduli $\bbC_\mathrm{KV}$)
and on the damper (with the viscous moduli tensor $\bbD$),
cf.\ the 5th row in Table~\ref{tab:models}. The {\it stored energy} is then
\begin{align}
E(e_\mathrm{el},\pi)
=\int_{\Omega}\!\Big(\frac12\bbC_\mathrm{M}e_\mathrm{el}{:}e_\mathrm{el}
+\frac12\bbC_\mathrm{KV}\pi{:}\pi\Big)\dd x
\end{align}
while the dissipation rate is $\bbD\DT\pi{:}\DT\pi$.
Abbreviating $\calE(u,\pi)=E(e(u){-}\pi,\pi)$, testing \eqref{difsystem-a}
by $\DT u$ and using the rheological ansatz \eqref{eq:GDE}
and the boundary conditions (\ref{difsystem}b,c), after a little calculus
one obtains the {\it total energy balance} in the form:
\begin{align}\label{TotalEngr-gen}
\calE(u(t),\pi(t))
+\int_0^t\!\!\int_{\Omega}\bbD\DT\pi{:}\DT\pi
\,\d x\d t
=\calE(u_0,\pi_0)
+\int_0^t\!\bigg(\int_{\Omega}\!f{\cdot}\DT u\,\d x+
\int_{\GNeu}\!\!g{\cdot}\DT u\dd S\bigg)\dd t.
\end{align}
For simplicity, here we assumed homogeneous Dirichlet condition $w=0$.
The time integrals on the left- and right-hand side of \eqref{TotalEngr-gen}
represent the dissipated energy due to viscosity and the work of external
forces done over the time interval $[0,t]$, respectively. Note that we need
to prescribe the initial conditions both $u(0,\cdot)=u_0$ and
$\pi(0,\cdot)=\pi_0$.
In a general case if $w\ne0$, one can first make a substitution
of $u{-}\bar w$ with an extension $\bar w$ of the boundary data $w$ inside
the bulk domain and then formulate an energy balance for a ``shifted''
solution satisfying homogeneous Dirichlet condition but with a modified
loading $f$ and $g$ while the internal variable $\pi$ remains unaffected.

As a special case, we can get both the Kelvin-Voigt model and
the Maxwell model. The former model results as the limit for
$\bbC_\mathrm{M}\to\infty$, which yields $e_\mathrm{el}=0$
so that simply $e(u)=\pi$ and, for $\bbD=\chi\bbC$, the energy balance
\eq{TotalEngr-gen} simplifies as
\begin{align}\label{TotalEngr}
\calE(u(t))+\int_0^t\!\!\int_{\Omega}\chi\bbC e(\DT u){:}e(\DT u)\,\d x\d t
=\calE(u_0)+\int_0^t\!\bigg(\int_{\Omega}\!f{\cdot}\DT u\,\d x+
\int_{\GNeu}\!\!g{\cdot}\DT u\dd S\bigg)\dd t,
\end{align}
with $\calE(u)=\int_{\Omega}\!\frac12\bbC e(u){:}e(u)\dd x$.
The Maxwell model results as the limit for $\bbC_\mathrm{KV}\to0$;
the splitting \eqref{e(u)=e-pi} and in particular the internal
variable $\pi$ remains in this model.

The other, higher-order models need more involved considerations and we will not
present it here. In particular, the 4-parameter solid uses again
\eqref{e(u)=e-pi} but the Burgers rheology, having two ``free nodes'' (cf.\
the rheological scheme at the 7th row in Table~\ref{tab:models}), needs
introduction of two internal variable and decomposition of $e(u)$
in \eqref{e(u)=e-pi} into 3 terms.

Under appropriate qualification of the external
loading and the initial conditions, energy balance \eq{TotalEngr}
gives also a-priori estimates of the solutions in respective norms
by using typically the Gronwall, the Young, and the H\"older inequalities.
Due to convexity of the energy $\calE(\cdot)$, this manipulation
can be reflected to the implicit time-discretisation schemes
considered in this paper, yielding numerical stability
and convergence of such schemes for $\tau\to0$. In our linear
situation, this convergence is indeed simple.

Evaluation and visualization of the spatial distribution of
the energies occurring in balances like \eqref{TotalEngr-gen} or \eqref{TotalEngr} may be of a
special interest,
since it shows in which regions of the body the dissipation takes place,
see the numerical example of Section~\ref{Simple} or \cite{RoPaMa_Visc}.
This energy dissipation leads to a heat production (not considered here,
however), and thus its spatial distribution would be important when solving the
heat-transfer problem in a possibly full thermomechanical coupling.
These forms of energetics are also of interest, since they could be used
to solve contact problems of visco-elastic bodies, see Section~\ref{Contact},
or even more complex problems where also inelastic phenomena take place on
the boundaries (or interfaces) of the viscous bodies, cf. \cite{RoPaMa_Visc}.
Techniques for the evaluation of these energies in combination with BEM have
been briefly described in Sections \ref{ViscousContact} and \ref{sect-other},
and employed in Section \ref{Numericalexamples}.

}

\end{document}